\documentclass[11pt]{amsart}
\usepackage{amssymb}
\usepackage{amsmath}
\usepackage{amsfonts}
\usepackage{pstricks}
\usepackage{pst-node}
\usepackage{graphicx}
\theoremstyle{plain}
\newtheorem{thm}{Theorem}[section]

\newtheorem{dfn}[thm]{Definition}
\newtheorem{rmk}[thm]{Remark}

\newtheorem{prb}[thm]{Problem}

\theoremstyle{remark}

\def\pmc#1{\setbox0=\hbox{#1}
    \kern-.1em\copy0\kern-\wd0
    \kern.1em\copy0\kern-\wd0}

\begin{document}

\bigskip

\title[On $F-$bubbles and hereditary $F-$functorial equivalence]{
On $F-$bubbles and hereditary $F-$functorial equivalence of topological spaces}

\bigskip

\author[U.~H.~Karimov]{Umed H. Karimov}
\address{Institute of Mathematics,
Academy of Sciences of Tajikistan, Ul. Ainy $299^A$,
Dushanbe 734063,
Tajikistan}
\email{umedkarimov@gmail.com}


\begin{abstract}
In \cite{KarRep1998}  the concept of $\check{H}^n-$bubles was defined and investigated. In this paper we generalize this conception for some other functors $F$. Open questions are formulated.
\end{abstract}

\date{\today}

\maketitle

\section{Introduction}

The $\check{H}^n$-bubbles ($\check{H}^*$ is \v{C}ech cohomology with integer coefficients) were defined and investigated in \cite{KarRep1998}. Topological space $X$ is $\check{H}^n$-bubble if $\check{H}^n(X) = \mathbf{Z}$ and $\check{H}^n(M) = \mathbf{0}$ for any proper subset $M \subset X.$
Every compact connected metrisable $n$-dimensional homology or cohomology manifold is $\check{H}^n$-bubble. If in $n$-dimensional sphere, $n > 1$, we identify any two point we get $\check{H}^n$-bubble which is not homology manifold.

We wish to generalize this concept for singular homology and other functors. According to the famous Barrat and Milnor example \cite{BM} there exists proper subset of $S^3$ with non-trivial $3-$dimensional singular homology group therefore the definition of $n-$dimensional bubble for singular homology is more complicated, needs the concept of hereditary equivalence of topological spaces. The purpose of this paper is to define the concept $H_n^s$-bubble for singular homology group, generalize it, and state some open problems.

\section{Main definitions and results}

Let $X$ and $Y$ be fixed  topological spaces and $F$ a covariant or contravariant topologically invariant functor (in sequel we will consider only topologically invariant functors) which is defined on all subspaces of these spaces with values in the category of abelian groups.

\begin{dfn} The spaces $X$ and $Y$ are called hereditary $F$-equivalent if:
\begin{enumerate}
  \item There exists bijection $f: X \rightarrow Y$,
  \item There exists a isomorphism $f_*(M): F(M) \rightarrow F(f(M))$ for every $M \subset X$ if $F$ is covariant functor and
        there exists a isomorphism $f^*(M): F(f(M)) \rightarrow F(M)$ if $F$ is contravariant functor,
  \item For every subsets $N\subset M\subset X$ the diagrams for covariant functors $F$\\
\[\begin{array}{lcr}
{F(M)}   & \stackrel{f_*(M)}\longrightarrow & {F(f(M))}   \\
\ \ \ \ \uparrow  &                          & \uparrow  \ \ \ \  \\
{F(N)}   & \stackrel{f_*(N)}\longrightarrow & {F(f(N))},    \\
\end{array}\]
and diagrams for contravariant functors $F$
\[ \begin{array}{lcr}
{F(M)}    & \stackrel{f^*(M)}\longleftarrow & {F(f(M))}    \\
\ \ \ \ \downarrow &                         & \downarrow \ \ \ \   \\
{F(N)}    & \stackrel{f^*(N)}\longleftarrow & {F(f(N))},    \\
\end{array} \]
in which vertical arrows are generated by inclusions are commutative.
\end{enumerate}
The mapping $f$ it is called hereditary equivalence of the spaces $X$ and $Y.$
\end{dfn}

\begin{dfn}
The space hereditary equivalent to $n$-dimensional sphere with respect to functor $F$  is called $n-$dimensional $F-$bubble.
\end{dfn}

In the case when the functor $F$ is the functor of $n-$dimensional \v{C}ech cohomology $\check{H}^n$ this definition obviously coincides with the definition of $\check{H}^n-$bubble $X$ ($f-$equivalence is any bijective mapping of the space $X$ and $n-$dimensional sphere $S^n).$

It is obvious that homeomorphic spaces are hereditary $F$-equivalent with respect to any topologically invariant functor $F$ which is defined on all subspaces of $X$ and $Y.$ Therefore non hereditary equivalent with respect some functor $F$ and for any bijection are not homeomorphic.

Let us show on example how this concept can be applied. Let $F$ be the functor of 1-dimensional singular homology with the coefficients in the group of integers and let $MB$ be the M\"{o}bius band, $C$ is cylinder $S^1\times [0; 1]$. The spaces $MB$ and $C$ are non $H_1^s-$equivalent with respect to any bijection because the compliment of any boundary point of $MB$ and only such points do not change the homology groups. Such points in $C$ are only boundary points, therefore for any hereditary equivalences $f$ for the functor $H_1^s$ the boundary points of $MB$ should correspond to boundary points of $C$ and vice versa, but the boundary $\partial (MB)$ of the M\"{o}bius band is homeomorphic to topological circle, $H_1^s(\partial (MB)) \approx \mathbf{Z} $ and boundary of cylindric surface $\partial C$ is the topological sum of two circles, $H_1^s(\partial C) \approx \mathbf{Z} \oplus \mathbf{Z},$ therefore these two spaces are not homeomorphic.

Subspaces of the real line $X = \{x\ |\ x = 0; 1; 2; \dots n; \dots\}_{n\in N}$ and $Y = \{x\ |\ x = 0; \frac{1}{1}; \frac{1}{2}; \dots \frac{1}{n}; \dots \}_{n\in N} $ with natural bijection are hereditary equivalent with respect to singular homology but these spaces are not homeomorphic.

\begin{dfn} The functor $F_1$ is more informative than the functor $F_2$ on the category of spaces $\mathcal{F}$ if for every hereditary $F_1-$equivalent spaces $X$ and $Y$ belonging to $\mathcal{F}$ it follows that the spaces $X$ and $Y$ are hereditary $F_2-$equivalent with the same bijection. The functor  $F_1$ is strongly more informative than the functor $F_2$ if $F_1$ is more informative than $F_2$ and there exist spaces $X_1$ and $Y_1$ in $\mathcal{F}$ which are  hereditary $F_2-$ equivalent but which are not hereditary  $F_1-$equivalent. The functors $F_1$ and $F_2$ are hereditary equivalent if $F_1$ is more informative than $F_2$ and vice versa. The functor $F$ is called complete functor if it is informative each other functor on the category $\mathcal{F}$.
\end{dfn}

For example on the category of zero-dimensional spaces the functor of zero-dimensional \v{C}ech cohomology $\check{H}^0$ is more informative than the functor of zero-dimensional singular cohomology $H_s^0$ since all $0-$dimensional metrisable spaces with the same cardinality have isomorphic $0-$dimensional singular homology groups (they are isomorphic to the direct product of the group $\mathbf{Z}$ in the quantity of the cardinality of the point of the space).The functor $\check{H}^0$ is strongly more informative than the functor $H_s^0$ since natural bijective mapping of the space $\{x|\ x = 0; 1; 2; \dots n; \dots\}_{n\in N} \subset \mathbf{R^1}$ to $\{x|\ x = 0; \frac{1}{1}; \frac{1}{2}; \dots \frac{1}{n} \dots \}_{n\in N}  \subset \mathbf{R^1} $ induces isomorphic mappings of the singular homologies of all subspaces but the $0$-dimensional \v{C}ech cohomology of these spaces are different, they are isomorphic to direct product and direct sum of the countable numbers of  $\mathbf{Z},$ respectively. Therefore the functor $H^s_0$ is not complete on the category of metrisable spaces but the following  theorem holds:

\begin{thm}\label{Sing complete}
On the category of finite polyhedra the functor of zero-dimensional singular homology $H^s_0(\ )$ is a complete functor.
\end{thm}
\begin{proof} Let $X$ and $Y$ be two finite polyhedra which are hereditary equivalent with respect to the functor  $H^s_0.$ Then there exists bijective mapping $f: X \ \rightarrow Y$ such that for every $M \subset X$ isomorphism $H^s_0(M)$ to $H^s_0(f(M))$ is defined. Let us prove that f is a homeomorphism.

Let $U$ be open in $Y$ and $f^{-1}(U)$ be not opened in $X$, i.e. there exist point $x \in f^{-1}(U)$ and sequence $x_n \in X\setminus f^{-1}(U)$ convergent to $x$. Without loss of generality we may assume that all $x_n$ belong to one simplex and if we connect by segments $x_1$ and $x_2,$ $x_2$ and $x_3$ and so on, and add the point $x$ to this set, we get arc $L.$ Since $Z = H_0^s(L) = H_0^s(f(L))$ the subspace $f(L)$ is linearly connected subspace of $Y$ and therefore is arcwise connected space (see e.g. \cite[\S 50, II, Theorem 1 and I, Theorem 2]{K}) and there exist injective continuous mapping $h$ of the segment $[0; 1]$ to $f(L)$ such that $h(0) = f(x_1), h(1) = f(x).$ The set $h^{-1}(U)$ is open in $[0; 1].$ Therefore there exists point $a \in [0; 1)$ such that $h([a; 1]) \subset U\cap f(L)$ and $h(a) \neq f(x)$. Then $f^{-1}(h([a, 1]))$ is linearly connected non one-point subspaces of the arc $L$, contained $x$ but not contained $x_k$ for any $k.$ This contradicts to the supposition that the sequence $\{x_k\}_{k\in N}$ converges to $x.$ Therefore $f^{-1}(U)$ is open, $f$ is bijective continuous mapping of compactum, $f$ is homeomorphism and $H^s_0(\ )$ is a complete functor.
\end{proof}

\begin{rmk}
The condition of finiteness of polyhedra is essential as the following example shows: Let polyhedron $X$ be the countable bouquet of $S^1,$ let $Y$ be a Hawaiian earring and let $f$ be a natural bijection $X \rightarrow Y$. These two spaces are not homeomorphic but they are obviously hereditary $H_0^s$ equivalent  and $f$ is $H_0^s$ equivalence.
\end{rmk}

Similarly to the proof of Theorm \ref{Sing complete}, it is possible to prove the following theorem:

\begin{thm}
Let $Q$ be Hilbert cube and $f: Q \mapsto Q$ hereditary equivalence with respect to the functor $H^s_n$ for some $n\geq 0$. Then $f$ is homeomorphism.
\end{thm}

\begin{rmk}
Any two full functors are hereditary equivalent and the relation of hereditary equivalence is reflexive, symmetric and transitive.
\end{rmk}

\begin{thm}
On the category of Hausdorff spaces with the first axiom of countability the group of continuous functions with the value in integers $C(\ ; \mathbf{Z})$ are full functor.
\end{thm}

\begin{proof}
Let $X$ and $Y$ be two Hausdorff spaces with the first axiom of countability which are hereditary equivalent with respect to the functor $C(\ ; \mathbf{Z})$. Let us prove that in this case the spaces $X$ and $Y$ are homeomorphic and therefore these functors are full functors.

Suppose that $f: X \rightarrow Y$ is hereditary equivalence. Let the sequence of points $x_n \in X$ converges to the point $a.$ Let $M = \{a, x_1, x_2, \dots, x_n, \dots\}_{n\in N}.$ Since $X$ is Hausdorff space the group $C(M ;\mathbf{Z})$ isomorphic to direct sum  $\sum_1^{\infty}\mathbf{Z}$ of countable number of the $\mathbf{Z}$. Let us prove that the sequence $f(x_n)$ converges to $f(a)$. Suppose that the sequence $\{f(x_n)\}_{n\in N}$ does not converge to $f(a)$ then without loss of generality, because $Y$ is Hausdorff spaces with the first axiom of countability, it is possible to assume that subspace $f(M)$ is discrete and the group $C(f(M); \mathbf{Z})$ isomorphic to the direct product $\prod_1^{\infty}\mathbf{Z}$. But  the groups $\sum_1^{\infty}\mathbf{Z}$ and $\prod_1^{\infty}\mathbf{Z}$ are not isomorphic.

\end{proof}

\begin{rmk}
The group $C(M ;\mathbf{Z})$ is isomorphic to the group of section of constant sheaf with the stalk $\mathbf{Z}$ which is natural isomorphic to $0$-dimensional sheaf cohomology and isomorphic to $0$-dimensional \v{C}ech cohomology \cite{B}. On the category of Hausdorff spaces with first axiom of countability subset $A\subset X$ is closed if and only if $\check{H}^0(M\cup x) = \check{H}^0(M) \oplus \mathbf{Z}$ for every $M \subset A$ and for every $x \notin A$ or in the terms of $1-$dimensional \v{C}ech cohomology, $A$ is closed in $X$ if and only if $\check{H}^1(M\cup x, M) = \mathbf{0}$ for every subset $M\subset A$ and every $x\notin A.$
\end{rmk}

\begin{thm}\label{S^1}
Metrisable space which are hereditary equivalent to the circle $S^1$ with respect to $1-$dimensional singular homology functor $H_1^s$  are homeomorphic to circle.
\end{thm}

\begin{proof} Let $f$ be hereditary equivalence of the circle $S^1$ and the space $X$ with respect to the functor $H_1^s.$ Then $H_s^1(X) = \mathbf{Z}$ and there exists continuous mapping $g: S^1 \rightarrow X,$ which induces isomorphism $g_1: H_1^s(S^1) \mapsto H_1^s(X).$ Consider the composition of the mappings
$$H_1^s(S^1) \rightarrow H_1^s(g(S^1)) \rightarrow H_1^s(X).$$
This composition is not trivial, therefore $H_1^s(g(S^1)) \neq 0.$ Since all proper subsets of the circle have trivial $1-$dimensional singular homology it follows that  $g(S^1) = X$ and the space $X$ as image of Peano continuum is Peano continuum. The space $X$ does not contains more then one simple closed curve since all proper subsets of $X$ $1-$acyclic, and contains simple closed curve since otherwise it should be dendrite and $AR$ space \cite[\S\ 53, Theorem 16]{K} but $H_1^s(X) \neq 0.$ Therefore $X$ contains exactly one simple closed curve and since all proper subsets of $X$ are acyclic it follows that $X$ itself is circle.
\end{proof}

\bigskip

\begin{thm}
Let $f: X \rightarrow Y$ be hereditary equivalence with respect to the functor $H_n^s,\ n \geq 1,$  of two finite polyhedra and every point $x \in X$ belongs to some topological simplex of the dimension $n + 1$ i.e. to the space homeomorphic the standard simplex of the dimension $n + 1$. Then $f$ is homeomorphism.
\end{thm}

\begin{proof}
Suppose that the sequence $\{x_k\}_{k\in \mathbf{N}}$ converges to the point $a$ but the sequence $\{f(x_k)\}_{k\in \mathbf{N}}$ does not converge to the point $f(a).$ Since $X$ is finite polyhedron then without loss of generality it is possible to suppose that all points $\{x_n\}_{n\in \mathbf{N}}$ are different, belong one $n+1-$dimensional simplex and it is possible to connect them by topological segments such that at the end we get arc $l_1$ which connects $x_1$ and $a$ and which contains all the points $\{x_k\}_{k\in \mathbf{N}}.$ In the $n+1-$dimensional simplex consider $n-$dimensional topological spheres $S^n_k$ each of which contains the corresponding point $x_k,$ which does not intersects each other and with the arc $l_1,$ diameters of which tends to $0.$ Obviously the subspace $l_1\cup (\cup_{k=1}^{\infty} S^n_k) $ has homotopical type of $n-$dimensional Hawaiian Earring $\mathcal{HE}^n$ and its $n-$dimensional singular homology group is uncountable. Consider the element $e\in H_n^s(l_1\cup (\cup_{k=1}^{\infty} S^n_k))$ the restriction of which in homology $H_n^s(S^n_k)$ is non trivial for every sphere $S^n_k.$ Since the polyhedron $Y$ is compact space it is possible to assume that the sequence $f(x_n)$ converges to some point $f(b), f(b) \neq f(a).$ Moreover it is possible to assume that $f(b)\notin f(l_1\cup (\cup_{k=1}^{\infty} S^n_k))$ because in any neighborhood of the point $a$ there exists subspace $l_m\cup (\cup_{k=m}^{\infty} S^n_k)$ for large enough $m$ homeomorphic to $l_1\cup (\cup_{k=1}^{\infty} S^n_k)$ ($l_m$ is part of $l_1$ connecting $x_m$ and $a$). The element $f(e)\in H_n^s(f(l_1\cup (\cup_{k=1}^{\infty} S^n_k)))$ as every element of the singular homology group has compact support $M.$ Compact subspace $M$ can contain only finite members of the sequence $f(x_n)$ because $f(b)\notin f(l_1\cup (\cup_{k=1}^{\infty} S^n_k)).$ Therefore $f^{-1}(M)$ can contain only finite number of spheres $S^n_k$ and the element $e$ does not belongs to the image of the mapping $H_n^s(f^{-1}(M)) \rightarrow H_n^s(l_1\cup (\cup_{k=1}^{\infty} S^n_k)).$ We get contradiction i.e. $f$ is bijective continuous mapping of compacta and $f$ is homeomorphism.
\end{proof}

In the paper \cite{KarRep1998} following theorem was proved:

\begin{thm}
Every compact metrisable $\check{H}^n$-bubble is Peano continuum.
\end{thm}

The similar statement is valid for singular $n-$bubles:

\begin{thm}
Every metrisable $H_n^s$-bubble $X$ is Peano continuum.
\end{thm}

\begin{proof} Obviously $H_n^s$-bubble $X$ is linearly connected space. Consider the generator of the group $1\in \mathbf{Z} \approx H_n^s(X)$ and the cycle which corresponds to this generator. The cycle is the sum of finite number of $n-$dimensional singular simplexes. The image of singular simplex is continuum Peano. Let $P^n$ is the union of these continua. Since the composition
$$H_n^s(P^n) \rightarrow H_n^s(f(P^n)) \rightarrow H_n^s(X)$$
is surjective mapping and every proper subset of $X$ is acyclic in the dimension $n$ then $f(P^n) = X.$ The union of finite number of Peano continua is locally connected and since the space $X$ is connected it follows that $X$ is Peano continuum.
\end{proof}

\section{Some open problems}

Theorem \ref{S^1} could be considered as homology characterization of the circle $S^1$. It is interesting to find the homology characterization of $S^2:$ Zastrow proved that every subset of the plane is acyclic with respect to $H_2^s$ therefore the usual $2-$sphere $S^2$ is $H_2^s$-bubble \cite{Z1997,Z}.

\begin{prb}
Is every $H_2^s$-bubble $X$ with trivial $1-$dimensional singular homology group homeomorphic to $S^2?$
\end{prb}

\begin{prb}
Does there exist $n$-dimensional, $ n > 2,$ or infinite dimensional continuum Peano
 $X$ all subsets $M$ of which:

\begin{itemize}
  \item Have trivial  homotopy groups $\pi_m(M)$ for some $m \geq n ?$
  \item Have trivial  homotopy groups $\pi_m(M)$ for all  $m \geq n ?$
\end{itemize}

\end{prb}

In \cite{KarRep1998} were proved two theorems the following wordings of which clarify and correct the corresponding formulations of the theorems of the paper \cite{KarRep1998}:

\begin{thm}\label{bubble}
There exist $2$-dimensional $2-$cyclic compact metrisable absolute neighborhood retract $ANR$ which does not contain any $\check{H}^2$-bubles.
\end{thm}

\begin{thm}
Every $n-$cyclic finite-dimensional locally $L\check{H}^n-$trivial metrisable compactum contains an $\check{H}^n-$bubble {\rm(}a topological space is said hereditarily locally cohomologically $n-$trivial {\rm(}$L\check{H}^n-$trivial{\rm)} if every point has neighborhood such that cohomology $\check{H}^n$ of all subspaces are trivial{\rm)}.
\end{thm}

\begin{prb}
Are a similar statements are valid for bubbles with respect to singular (co)-homology ?
\end{prb}

Theorem \ref{bubble} is interesting to compare with Theorem 5.5 of \cite{Kup} which states that any $q-$dimensional, $q-$cyclic compact $ANR$ contains a $q-$bubble.

\section{Acknowledgements}

I am grateful to Professor G. Conner for the invitation to the Arches Topology Conference (https://math.byu.edu/~curtkent/arches/2018.html) and for the discussion during this event of the main result of the paper. I am grateful to D. Repov\v{s} for the remarks and comments.

\end{document}